\documentclass{amsart}
\usepackage{amsfonts}
\usepackage{amsmath}
\usepackage{amssymb}
\usepackage{amsthm}

\newtheorem{lemma}{Lemma}
\newtheorem{theorem}{Theorem}
\newtheorem{rem}{Remark}
\newtheorem{exmp}{Example}

\begin{document}

\title{On Tits buildings of type $A_{n}$}
\author{Mark  Pankov}
\address{Institute of Mathematics NASU, Kiev}
\email{pankov@imath.kiev.ua} 

\begin{abstract}
Let $P$ and $P'$ be projective spaces having the same dimension,
this dimension is denoted by $n$ assumed to be finite.
Denote by ${\mathfrak F}$ and ${\mathfrak F}'$ the sets
of maximal flags of the spaces $P$ and $P'$, respectively.
A subset of ${\mathfrak F}$ (${\mathfrak F}'$) is said to be
an {\it apartment} if it is the intersection of ${\mathfrak F}$ (${\mathfrak F}'$)
and an apartment of the $A_{n}$-building associated with the projective space $P$ ($P'$).
We show that any mapping $f:{\mathfrak F}\to {\mathfrak F}'$
sending apartments to apartments is induced by a strong embedding of $P$ to $P'$ or
to the dual space $P'^{*}$.
Moreover this embedding is a collineation if $f$ is surjective.
\end{abstract}

I dedicate this paper to my wife Inna and my daugther 
Sabina-Stefany.

\maketitle

\section{Introduction}
According to J. Tits \cite{Tits1}, a building is a chamber complex
together with a distinguished family of subcomplexes called
{\it apartments}\, and satisfying certain collection of axioms.
In the present paper we will consider so-called $A_{n}$-buildings
(buildings of type $A_{n}$); these buildings are associated with projective spaces.
We study mappings of the chamber sets of $A_{n}$-buildings
which send apartments to apartments.
The main result of the paper (Theorem 3.1) says that
these mappings are induced by strong embeddings of the corresponding projective spaces.

\subsection{Chamber complexes}
Let $X$ be a set and $\Delta$ be a set of proper subsets of $X$.
We say that $\Delta$ is a {\it complex} over $X$ if
the following conditions hold true:
\begin{enumerate}
\item[(C1)] $\Delta$ contains each one-element subset of $X$,
\item[(C2)] if $A\in \Delta$ then any proper subset of $A$ belongs to $\Delta$.
\end{enumerate}
For this case elements of $\Delta$ are known as {\it simplexes},
one-element simplexes (elements of the set $X$) are said to be {\it vertices}
of the complex $\Delta$.

Recall that the {\it rank} of a complex
is the maximum of the cardinal numbers of simplexes belonging to this complex.
{\it In this paper we will always require that the rank of a complex is finite.}

Let $\Delta$ be a complex over a set $X$,
let also $X'$ be a subset of $X$ and $\Delta'$ be a complex over $X'$.
If $\Delta'$ is contained in $\Delta$ then we say that
$\Delta'$ is a {\it subcomplex} of the complex $\Delta$.

Two complexes $\Delta$ and $\Delta'$ over sets $X$ and $X'$ (respectively)
are said to be {\it isomorphic} if there exists a bijection $f:X\to X'$
such that $f(\Delta)=\Delta'$,
this bijection is called an {\it isomorphism}.

\begin{exmp}[Flag complex]{\rm
Let $(X,\le)$ be a partially ordered set.
A set $F\subset X$ is said to a {\it flag} if for any two elements
$x,y\in F$ we have $x\le y$ or $y \le x$.
The set of all flags is a complex over $X$, it is known as
the {\it flag complex} of the partially ordered set $(X,\le)$.
}\end{exmp}

Now consider a complex $\Delta$ which rank is finite and denoted by $n$.
Then each simplex is contained in a maximal simplex;
maximal simplexes are called {\it chambers}.
We will also require that our complex satisfies the following condition:
{\it any two chambers have the same cardinal number}
(clearly, this number is equal to $n$).
Simplexes spanned by $n-1$ vertices are known as {\it panels}.
Two chambers are said to be  {\it adjacent} if their intersection is a panel.
We say that $\Delta$ is a {\it chamber complex} if
for any two chambers $C$ and $C'$ there is a sequence of chambers
$$C=C_{0},C_{1},\dots,C_{m}=C'$$
such that
$C_{i-1}$ and $C_{i}$ are adjacent for each number $i=1,\dots, m$.

A chamber complex is called {\it thick}
if each panel is contained in at least three chambers.
For the case when any panel is contained in exact two chambers
we say that our complex is {\it thin}.

\begin{exmp}[$A_{n}$-complex]{\rm
Set $I_{n}:=\{1,\dots,n+1\}$ and consider
the set of all proper subsets of $I_{n}$ ordered by the relation $\subset$.
The corresponding flag complex is denoted by $A_{n}$.
It is a thin chamber complex containing $(n+1)!$ distinct chambers.
The rank of this complex is equal to $n$.
}\end{exmp}

\subsection{Buildings}
Let $\Delta$ be a complex together with a family of subcomplexes
$\aleph$ called {\it apartments} and satisfying the following axioms:
\begin{enumerate}
\item[(TB1)]
All apartments are isomorphic thin chamber complexes.
\item[(TB2)]
For any two simplexes there exists an apartment containing them.
\end{enumerate}
Then $\Delta$ is a chamber complex
which rank is equal to the rank of the apartments.
The pair $(\Delta, \aleph)$ is said to be a {\it building} \cite{Tits1}
if the chamber complex $\Delta$ is thick and the following axiom holds true:
\begin{enumerate}
\item[(TB3)]
Let $\Sigma$ and $\Sigma'$ be apartments;
let also $S$ and $U$ be simplexes belonging to both $\Sigma$ and $\Sigma'$.
Then there exists an isomorphism of the complex $\Sigma$ to the complex $\Sigma'$
which leaves fixed all simplexes contained in $S$ and $U$.
\end{enumerate}

All buildings of finite rank $\ge 3$ with finite apartments
were completely determined by J. Tits \cite{Tits2}.
A lot of information related with buildings and their applications
can be found in \cite{Buekenhout}, see also \cite{Garret} and \cite{Taylor}.

In this paper we will consider only $A_{n}$-buildings, $n\ge 2$.
By the definition, a building $(\Delta, \aleph)$ has type $A_{n}$
if its apartments are isomorphic to the complex $A_{n}$ (Example 1.2).
It is well-known that for this case
$\Delta$ is the flag complex of certain $n$-dimensional projective space
\cite{Tits1}.

\subsection{Main result}
Let $(\Delta, \aleph)$ and $(\Delta', \aleph')$ be buildings of type $A_{n}$.
Let also $P$ and $P'$ be the associated projective spaces
(i.e. $\Delta$ and $\Delta'$ are the flag complexes of these spaces).
Put ${\mathfrak F}$ and ${\mathfrak F}'$ for the sets of all chambers of
the complexes $\Delta$ and $\Delta'$, respectively
(it is trivial that ${\mathfrak F}$ and ${\mathfrak F}'$ consist
of all maximal flags of $P$ and $P'$).
We consider mappings of ${\mathfrak F}$ to ${\mathfrak F}'$
which send apartments to apartments
(a subset of ${\mathfrak F}$ or ${\mathfrak F}'$ is said to be an {\it apartment}
if it consists of all chambers belonging to certain apartment of the corresponding building)
and show that these mappings are induced by
strong embeddings of $P$ to $P'$ or to the dual projective space $P'^{*}$ (Theorem 3.1).

\begin{rem}{\rm
Adjacency preserving bijections of ${\mathfrak F}$ to itself
were studied by H. Hav\-licek, K. List and C. Zanella \cite{HLZ}.
The background for the investigation of adjacency preserving mappings can be found in
\cite{Benz} and \cite{Lester}.
}\end{rem}

\section{Preliminaries: Basic properties of projective spaces}
\setcounter{exmp}{0}
\setcounter{rem}{0}

\subsection{Projective spaces}
A {\it projective space} is a set of points $P$
together with a family of proper subsets called {\it lines}
and satisfying certain collection of axioms (see, for example, \cite{VeblenYoung}).
One of these axioms says that for any two distinct points
$p,q\in P$ there is unique line containing them, we will denote this line by $p\,q$.

A set $S\subset P$ is said to be a {\it subspace}
if for any two distinct points $p$ and $q$ belonging to $S$
the line $p\,q$ is contained in $S$.
By this definition, the empty set and one-point sets are subspaces
(since they do not contain two distinct points).
The intersection of any collection of subspaces is a subspace.
For any subset $X\subset P$ the minimal subspace containing $X$
(the intersection of all subspaces containing $X$)
is called {\it spanned} by $X$ and denoted by $\overline{X}$.

We say that a subspace is $n$-{\it dimensional} if
$n+1$ is the smallest number of points spanning this subspace;
a subspace is {\it infinity-dimensional} if it is not spanned by
a finite set of points.
The empty set is unique $(-1)$-dimensional subspace,
points and lines are $0$-dimensional and $1$-dimensional subspaces,
respectively. The dimensions of other subspaces are greater than $1$.
Note also that the dimension of a projective space is not less than $2$
(otherwise, there is only one line which coincides with the space).

A set $X\subset P$ is said to be {\it independent} if
the subspace $\overline{X}$ is not spanned by a proper subset of $X$.
An independent set $B$ contained in a subspace $S$ is called a {\it base} for
$S$ if $\overline{B}=S$;
for the case when $S$ coincides with $P$ we will say that
$B$ is a base for our projective space.

It is well-known that for any subspace $S\subset P$ the following statements hold true:
\begin{enumerate}
\item[{\rm (1)}]
If $S$ is $n$-dimensional then any base for $S$ consists of $n+1$ points.
\item[{\rm (2)}]
Each base for $S$ can be extended to a base for any subspace containing $S$.
\end{enumerate}
By (2), for any two subspaces $S$ and $U$
the inclusion $S\subset U$ implies that
the dimension of $S$ is not greater than the dimension of $U$
and these subspaces have the same dimension only for the case when $S=U$.

\begin{exmp}{\rm
Let $P$ be an $n$-dimensional projective space, $n\ge 3$
and ${\mathcal L}$ be the set of lines of this space.
Let also $p\in P$.
Denote by ${\mathcal L}_{p}$ the set of all lines containing the point $p$.
A subset of ${\mathcal L}_{p}$ is said to be a {\it line}
if it consists of all lines $L\in {\mathcal L}_{p}$ contained
in certain plane ($2$-dimensional subspace) of $P$.
The set ${\mathcal L}_{p}$ together with the family of lines
defined above is an $(n-1)$-dimensional projective space.
A subset of ${\mathcal L}_{p}$ is a $k$-dimensional subspace
if and only it consists of all lines $L\in {\mathcal L}_{p}$
contained in certain $(k+1)$-dimensional subspace of $P$.
}\end{exmp}

\begin{exmp}{\rm
Let $V$ be a left vector space over a division ring.
Denote by $P(V)$ the set of all $1$-dimensional subspaces of $V$.
A subset of $P(V)$ is a {\it line} if it consists of all
$1$-dimensional subspaces contained in certain
$2$-dimensional subspace of $V$.
The set $P(V)$ together with the family of lines defined above
is a projective space if the dimension of $V$ is not less than $3$.
This projective space is $n$-dimensional if the dimension of $V$ is equal to $n+1$.
}\end{exmp}

Let $P$ and $P'$ be projective spaces,
let also ${\mathcal L}$ and ${\mathcal L}'$ be the sets of lines.
These spaces are called {\it isomorphic}
if there exists a bijection $f:P\to P'$ such that
$$f({\mathcal L})={\mathcal L}';$$
this bijection is said to be a {\it collineation}.

\begin{rem}{\rm
Projective spaces isomorphic to the projective spaces associated with vector spaces
are known as {\it Desarguesian}.
Any projective space the dimension of which is not less than $3$
is Desarguesian (see \cite{Baer} or \cite{VeblenYoung}).
However, non-Desarguesian planes ($2$-dimensional projective spaces) exist.
}\end{rem}

\subsection{Dual space}
Let $P$ be a projective space.
The dimension of this space is assumed to be finite and denoted by $n$.
Put $P^{*}$ for the set of all $(n-1)$-dimensional subspaces of $P$.
We say that a subset of $P^{*}$ is a {\it line} if it consists of
all $(n-1)$-dimensional subspaces containing
certain $(n-2)$-dimensional subspace of $P$.
The set $P^{*}$ and the family of all lines defined above form an
$n$-dimensional projective space;
it is called {\it dual} to the projective space $P$.
We will exploit the following well-known properties of the dual space:
\begin{enumerate}
\item[(1)]
A subset of $P^{*}$ is a $k$-dimensional subspace if and only if
it consists of all $(n-1)$-dimensional subspaces containing
certain $(n-k-1)$-dimensional subspace of $P$.
Hence there is a one-to-one correspondence
between the set of all $k$-dimensional subspaces of $P^{*}$ and
the set of all $(n-k-1)$-dimensional subspaces of $P$.
In particular, the second dual space $P^{**}$ coincides with $P$.
\item[(2)]
Points $p^{*}_{1},\dots, p^{*}_{n+1}\in P^{*}$ form a base for
the dual space if and only if there is a base for $P$
such that each $p^{*}_{i}$ is spanned by points of this base.
\end{enumerate}

\begin{exmp}{\rm
Let $n\ge 3$ and $p\in P$.
The set ${\mathcal L}_{p}$ (consisting of all lines of $P$ passing through the point $p$)
has the natural structure of an $(n-1)$-dimensional projective space (Example 2.1).
The dual projective space consists of all
$p^{*}\in P^{*}$ containing $p$;
it is an $(n-1)$-dimensional subspace of $P^{*}$.
We will denote this subspace by $P^{*}_{p}$.
}\end{exmp}

\begin{exmp}{\rm
Let $V$ be a left vector space over a division ring,
the dimension of $V$ is assumed to be finite and not less than $3$.
Then $P(V)^{*}$ (the projective space dual to $P(V)$)
is isomorphic to $P(V^{*})$
($V^{*}$ is the vector space dual to $V$) \cite{Baer}.
}\end{exmp}

\subsection{Embeddings of projective spaces}
Let us consider projective spaces $P$ and $P'$
and denote by ${\mathcal L}$ and ${\mathcal L}'$ the sets of lines of these spaces.
We say that an injection $f:P\to P'$ is an {\it embedding}
if it is collinearity preserving and
maps each triple of non-collinear points
to non-collinear points, in other words, if the following two conditions hold true:
\begin{enumerate}
\item[(1)] for any line $L\in {\mathcal L}$ there exists
a line $L'\in {\mathcal L}'$ such that $f(L)\subset L'$,
\item[(2)]
for each line $L'\in {\mathcal L}'$
there is at most one line $L\in {\mathcal L}$
such that $f(L)\subset L'$.
\end{enumerate}
It is trivial that any bijective embedding is a collineation.

An embedding is said to be {\it strong} if it transfers
independent subsets to independent subsets.

\begin{rem}{\rm
Non-strong embeddings exist
(see, for example, \cite{BrezuleanuRadulescu} or \cite{Kreuzer2}).
}\end{rem}

\begin{rem}{\rm
Assume that our projective spaces have the same dimension.
Then any strong embedding of $P$ to $P'$ sends bases to bases.
Inversely, W.-l. Huang and A. Kreuzer \cite{HuangKreuzer} (see also \cite{Kreuzer1})
have shown that any base preserving surjection of $P$ to $P'$
is a collineation.
It is natural to ask are non-surjective base preserving mappings of $P$ to $P'$
strong embeddings?
}\end{rem}

Let $f:P\to P'$ be a strong embedding.
Then for any subspace $S\subset P$ the dimension of the subspace
$\overline{f(S)}$ is equal to the dimension of $S$.
Assume also that our projective spaces have the same dimension which is finite.
Then the embedding $f$ induces the mapping
$$f^{*}:P^{*}\to P'^{*}$$
which transfers $p^{*}\in P^{*}$ to $\overline{f(p^{*})}\in  P'^{*}$.
This mapping is a strong embedding of $P^{*}$ to $P'^{*}$;
it is said to be {\it dual} to $f$.
A direct verification shows that the second dual embedding coincides with
$f$, in other words, $f^{**}=f$.

\begin{exmp}{\rm
If the dimension of our spaces is not less than $3$ then
for any point $p\in P$ the embedding $f$ defines the strong embedding of
${\mathcal L}_{p}$ to ${\mathcal L}'_{f(p)}$ such that
the dual embedding is the restriction of $f^{*}$ to $P^{*}_{p}$.
}\end{exmp}

\begin{exmp}{\rm
Let $V$ and $V'$ be left vector spaces over division rings $R$ and $R'$ (respectively),
the dimensions of $V$ and $V'$ are assumed to be not less than $3$.
Then $P(V)$ and $P(V')$ are projective spaces (Example 2.2).
A mapping $l:V\to V'$ is called {\it semilinear} if
$$l(x+y)=l(x)+l(y)\;\;\;\;\;\forall\;x,y\in V$$
and there exists a homomorphism $\sigma: R\to R'$ such that
$$l(ax)=\sigma(a)l(x)\;\;\;\;\;\forall\;x\in V,a\in F;$$
it is easy to see that $\sigma$ is a monomorphism if the mapping $l$ is non-zero.
Any semilinear injection $l:V\to V'$ induces the mapping
$$P(l):P(V)\to P(V')$$
which sends each $1$-dimensional subspace $Rx$, $x\in V-\{0\}$ to
the subspace $R'l(x)$.
Note that the mapping $P(l)$ is non-injective, for example, if
$l$ is the semilinear injection of ${\mathbb R}^{2n}$ to ${\mathbb C}^{n}$
defined by the formula
$$l(x_{1},y_{1},\dots, x_{n},y_{n}):=(x_{1}+iy_{1},\dots, x_{n}+iy_{n}).$$
However, if $l$ preserves the linear independence
(maps any set of linearly independent vectors to a set of linearly independent vectors)
then $P(l)$ is a strong embedding.
Inversely,
{\it any strong embedding of $P(V)$ to $P(V')$ is
induced by a semilinear injection of $V$ to $V'$ preserving the linear independence};
it is a simple consequence of Faure-Fr\"{o}lisher-Havlicek's version of
the Fundamental Theorem of Projective Geometry
\cite{FaureFrolisher}, \cite{Faure} and \cite{Havlicek1}.
The classical version of this theorem (see \cite{Artin}, \cite{Baer}) says that
any collineation of $P(V)$ to $P(V')$ is
induced by a semilinear isomorphism of $V$ to $V'$;
this statement was first proved by O. Veblen \cite{Veblen}
for projective spaces over finite fields
(more historical information can be found in \cite{KarzelKroll}).
}\end{exmp}

\section{$A_{n}$-buildings and their morphisms}
\setcounter{exmp}{0}
\setcounter{rem}{0}

\subsection{Buildings associated with projective spaces}
Let $P$ be an $n$-dimensi\-onal projective space.
Denote by $\Delta$ the flag complex of $P$
(the flag complex of the set of all proper subspaces of $P$ ordered by $\subset$).
Let $B$ be a base for $P$ and
${\mathcal B}$ be the set of all subspaces spanned by points of this base.
The flag complex of $({\mathcal B},\subset )$ is isomorphic to $A_{n}$
and called the {\it apartment} associated with the base $B$.
{\it The complex $\Delta$ together with the family of all apartments
defined above is a building of type $A_{n}$}.

It was noted above that there are not other non-trivial building of this type:
{\it each $A_{n}$-building $(n\ge 2)$ is associated with certain
$n$-dimensional projective space}
\cite{Tits1}.

Denote by ${\mathfrak F}$ the set of all maximal flags
(chambers of the complex $\Delta$).
Since the dimension of the projective space is equal to $n$,
any maximal flag consists of $n$ proper subspaces;
it is a chain of subspaces
$$p\in S_{1}\subset \dots \subset S_{n-2}\subset p^{*}$$
where $p\in P$, $p^{*}\in P^{*}$ and each $S_{i}$ is an $i$-dimensional subspace.
Recall that a subset ${\mathfrak A}\subset {\mathfrak F}$ is called an {\it apartment}
in ${\mathfrak F}$ if it is the intersection of ${\mathfrak F}$ with
an apartment of $\Delta$;
in other words, there is a base for our projective space
such that ${\mathfrak A}$ consists of all maximal flags
which components are spanned by points of this base.

Let ${\mathfrak X}\subset {\mathfrak F}$.
The set consisting of all subspaces contained in flags belonging to ${\mathfrak X}$
will be called the {\it trace} of ${\mathfrak X}$ and denoted by $T({\mathfrak X})$.
If $B$ is a base for $P$ then the apartment of ${\mathfrak F}$ associated with $B$
is the maximal subset of ${\mathfrak F}$ which trace
coincides with the set of all subspaces spanned by points of $B$.

\begin{exmp}{\rm
For a subspaces $S\subset P$ we denote by ${\mathfrak F}(S)$
the set of all maximal flags containing $S$.
Let $p\in P$.
If $n\ge 3$ then  ${\mathfrak F}(p)$ is the set of all maximal flags of
the projective space ${\mathcal L}_{p}$ (Example 2.1);
apartments of ${\mathfrak F}(p)$ are the intersections
of ${\mathfrak F}(p)$ with apartments ${\mathfrak A}\subset {\mathfrak F}$
such that $p\in T({\mathfrak A})$.
}\end{exmp}

\subsection{Morphisms of $A_{n}$-buildings}
Now assume that $P'$ is another $n$-dimensi\-onal projective space.
The set of all maximal flags of this space will be denoted by ${\mathfrak F}'$.

Let $f:P\to P'$ be a strong embedding.
Since for any subspace $S\subset P$ the dimension of the subspace
$\overline{f(S)}$ is equal to the dimension of $S$,
the embedding $f$ induces the mapping of ${\mathfrak F}$ to ${\mathfrak F}'$
which sends any maximal flag
$${\mathcal F}=(p,S_{1},\dots,S_{n-2},p^{*})$$
to the maximal flag
$$f({\mathcal F})=(f(p),\overline{f(S_{1})},\dots, \overline{f(S_{n-2})},f^{*}(p^{*}));$$
This mapping is injective and transfers apartments to apartments,
besides it is bijective only for the case when $f$ is a collineation.

The dual principles of Projective Geometry (Subsection 2.3)
show that {\it the building of any finite-dimensional projective space
coincides with the building associated with the dual projective space}.
Therefore, any strong embedding of $P$ to $P'^{*}$ induces certain
injection of ${\mathfrak F}$ to ${\mathfrak F}'$ which
maps apartments to apartments.

\begin{theorem}
Any mapping of ${\mathfrak F}$ to ${\mathfrak F}'$ sending apartments to apartments
is the injection induced by a strong embedding of $P$ to $P'$ or $P'^{*}$;
in particular,
any surjection of ${\mathfrak F}$ to ${\mathfrak F}'$
transferring apartments to apartments is induced
by a collineation of $P$ to $P'$ or $P'^{*}$.
\end{theorem}

\begin{rem}{\rm
It was established by H. Havlicek, K. List and C. Zanella \cite{HLZ}
that if $n=3$ and $f$ is a bijective transformation of
${\mathfrak F}$ preserving the adjacency relation in both directions
then $f$ is induced by a collineation of $P$ to itself or to the dual projevive space.
}\end{rem}

\begin{rem}{\rm
For each number $k=1,\dots n-2$ denote by
${\mathcal G}_{k}$ and ${\mathcal G}'_{k}$ the Grassmann spaces of
$k$-dimensional subspaces of our projective spaces.
If $B$ is a base for $P$ then the set
consisting of all $k$-dimensional subspaces spanned by points of $B$
is known as the {\it base} subset of ${\mathcal G}_{k}$ associated with $B$.
Any strong embedding of $P$ to $P'$
induces certain injective mapping of ${\mathcal G}_{k}$ to ${\mathcal G}'_{k}$.
Similarly, strong embeddings of $P$ to $P'^{*}$
induce injections of ${\mathcal G}_{k}$ to ${\mathcal G}'_{n-k-1}$;
these are injections to ${\mathcal G}'_{k}$ if $n=2k+1$.
It is not difficult to see that these mappings
send base subsets to base subsets.
M. Pankov \cite{Pankov3} (see also \cite{Pankov1} and \cite{Pankov2})
has shown that for the case when $n\ge 3$
there are not other mappings of ${\mathcal G}_{k}$ and ${\mathcal G}'_{k}$
satisfying this condition;
i.e. if $n\ge 3$ and $f:{\mathcal G}_{k}\to {\mathcal G}'_{k}$
maps base subsets to base subsets then
$f$ is induced by a strong embedding of $P$ to $P'$ or $P'^{*}$
(the second possibility can be realized only for the case when $n=2k+1$);
this embedding is a collineation if $f$ is surjective.
Geometrical transformations of Grassmann spaces were studied
by many authors, see, for example,
\cite{Benz}, \cite{BlunckHavlicek}, \cite{Chow},
\cite{Havlicek2}, \cite{Huang}, \cite{Kreuzer2}, \cite{Zanella}.
}\end{rem}

\section{Proof of Theorem 3.1}
\setcounter{lemma}{0}

Let $f:{\mathfrak F}\to {\mathfrak F}'$ be a mapping
which transfers apartments to apartments.
First of all note that {\it $f$ is injective}
(for any two elements of ${\mathfrak F}$ there exists an apartment containing them;
by our hypothesis, the restriction of $f$ to any apartment is injective).
We will show that $f$ is induced by a strong embedding of $P$ to $P'$ or $P'^{*}$.

\subsection{Main Lemma}
Let ${\mathfrak A}$ be an apartment in ${\mathfrak F}$.
Then ${\mathfrak A}':=f({\mathfrak A})$ is an apartment in ${\mathfrak F}'$.
Let also
$$B=\{p_{1},\dots,p_{n+1}\}\;\mbox{ and }\;B'=\{p'_{1},\dots,p'_{n+1}\}$$
be the bases associated with ${\mathfrak A}$ and ${\mathfrak A}'$, respectively.
We set
$$p^{*}_{i}:=\overline{B-\{p_{i}\}}\;\mbox{ and }\;p'^{*}_{i}:=\overline{B'-\{p'_{i}\}}$$
for each number $i\in I_{n}$
\footnote{recall that $I_{n}=\{1,\dots, n+1\}$, see Example 1.2}
and denote by ${\mathfrak A}_{i}$ and ${\mathfrak A}^{i}$
(${\mathfrak A}'_{i}$ and ${\mathfrak A}'^{i}$)
the sets of all flags belonging to ${\mathfrak A}$ (${\mathfrak A}'$)
and containing $p_{i}$ or $p^{*}_{i}$ ($p'_{i}$ or $p'^{*}_{i}$), respectively.
It is easy to see that these sets contain $n!$ elements.

Note also that ${\mathfrak A}$ can be presented as the disjoint union of
${\mathfrak A}_{1},\dots, {\mathfrak A}_{n}$ and
the disjoint union of ${\mathfrak A}^{1},\dots, {\mathfrak A}^{n}$.
The intersection of ${\mathfrak A}_{i}$ and ${\mathfrak A}^{j}$
will be denoted by ${\mathfrak A}^{j}_{i}$.
This set is not empty if $i\ne j$, for this case it contains $(n-1)!$ elements.
Each ${\mathfrak A}_{i}$ (${\mathfrak A}^{i}$)
is the disjoint union of all ${\mathfrak A}^{j}_{i}$
(${\mathfrak A}^{i}_{j}$) such that $j\ne i$.

\begin{lemma}
There exists a permutation $\sigma_{f}$ of the set $I_{n}$ such that
one of the following possibilities is realized:
\begin{enumerate}
\item[(1)] $f({\mathfrak A}_{i})={\mathfrak A}'_{\sigma(i)}$
and $f({\mathfrak A}^{i})={\mathfrak A}'^{\sigma(i)}$ for all $i\in I_{n}$,
\item[(2)] $f({\mathfrak A}_{i})={\mathfrak A}'^{\sigma(i)}$
and $f({\mathfrak A}^{i})={\mathfrak A}'_{\sigma(i)}$ for all $i\in I_{n}$.
\end{enumerate}
\end{lemma}

This statement will be proved in a few steps (Subsection 4.2 -- 4.5).

\subsection{Maximal inexact subsets of apartments}
We say that a subset ${\mathfrak X}$ of the apartment ${\mathfrak A}$ is {\it exact}
if ${\mathfrak A}$ is unique apartment containing ${\mathfrak X}$;
otherwise, ${\mathfrak X}$ is said to be {\it inexact}.
It is easy to see that ${\mathfrak X}$ is exact if and only if
for each $i\in I_{n}$ the intersection of all subspaces
$S\in T({\mathfrak X})$ containing $p_{i}$ coincides with $p_{i}$.

Take two distinct numbers $i,j\in I_{n}$
and denote by ${\mathcal X}_{ij}$ the set of all subspaces
$S\in T({\mathfrak A})$ such that
$$p_{i}p_{j}\subset S\;\mbox{ or }\;p_{i}\notin S.$$
Put ${\mathfrak X}_{ij}$ for the set of all flags
spanned by subspaces belonging to ${\mathcal X}_{ij}$;
in other words, ${\mathfrak X}_{ij}$ is the maximal subset of ${\mathfrak F}$
which trace is ${\mathcal X}_{ij}$.
If $k\ne i$ then the intersection of all subspaces belonging to
${\mathcal X}_{ij}=T({\mathfrak X}_{ij})$ and containing $p_{k}$
coincides with $p_{k}$.
Any flag ${\mathcal F}\in {\mathfrak A}-{\mathfrak X}_{ij}$
contains a subspace intersecting the line $p_{i}p_{j}$ by $p_{i}$.
Therefore,
${\mathfrak X}_{ij}\cup \{{\mathcal F}\}$ is an exact
subset of ${\mathfrak A}$ and
{\it ${\mathfrak X}_{ij}$ is a maximal inexact subset of ${\mathfrak A}$}.

Now we show that there are not other maximal inexact subsets of
apartments.

\begin{lemma}
If ${\mathfrak X}$ is a maximal inexact subset of ${\mathfrak A}$
then there exist $i,j\in I_{n}$
such that ${\mathfrak X}_{ij}$ coincides with ${\mathfrak X}$.
\end{lemma}

\begin{proof}
Since ${\mathfrak X}$ is inexact,
there is a number $i\in I_{n}$ such that the intersection of all subspaces
$S\in T({\mathfrak X})$ containing $p_{i}$ does not coincide with $p_{i}$;
we denote this intersection by $U$.
Then one of the following possibilities is realized:
\begin{enumerate}
\item[(1)] the dimension of $U$ is greater than $0$,
\item[(2)] $U=\emptyset$.
\end{enumerate}
For the case (1) the inclusion $p_{i}p_{j}\subset U$ holds for
certain number $j\in I_{n}-\{i\}$.
For the second case we can take any number $j\in I_{n}-\{i\}$.
It is trivial that for each of these cases
$T({\mathfrak X})$ is contained in ${\mathcal X}_{ij}$.
By the definition of the set ${\mathfrak X}_{ij}$,
we have ${\mathfrak X}\subset {\mathfrak X}_{ij}$.
These are both maximal inexact subsets and the inverse inclusion holds true.
\end{proof}

In particular, Lemma 4.2 implies that
{\it all maximal inexact subsets of apartments of ${\mathfrak F}$ and ${\mathfrak F}'$
have the same cardinal number.}

\begin{lemma}
The mapping $f$ transfers inexact subsets of ${\mathfrak A}$
to inexact subsets of ${\mathfrak A}'$;
moreover, maximal inexact subsets go over to maximal inexact subsets.
\end{lemma}

\begin{proof}
Let ${\mathfrak X}$ be an inexact subset of ${\mathfrak A}$.
Then ${\mathfrak X}$ is contained in at least two distinct apartments of ${\mathfrak F}$.
The $f$-images of these apartments are distinct apartments
of ${\mathfrak F}'$ containing $f({\mathfrak X})$ (since $f$ is injective).
Hence $f({\mathfrak X})$ is an inexact subset of ${\mathfrak A}'$.

Now assume that the inexact subset ${\mathfrak X}$ is maximal
and consider a maximal inexact subset ${\mathfrak Y}\subset{\mathfrak A}'$
containing $f({\mathfrak X})$.
Since $f$ is injective, the remark given after Lemma 4.2 shows that
$f({\mathfrak X})$ and ${\mathfrak Y}$
have the same cardinal number.
Then $f({\mathfrak X})$ coincides with ${\mathfrak Y}$.
\end{proof}

\subsection{Complements to maximal inexact subsets}
In this subsection we will study the structure of the {\it complement} set
$${\mathfrak C}_{ij}:={\mathfrak A}-{\mathfrak X}_{ij}.$$

Any permutation $\sigma$ of the set $I_{n}$ defines the flag
$${\mathcal F}_{\sigma}=(p_{\sigma(1)},S_{1},\dots, S_{n-2},p^{*}_{\sigma(n+1)})
\in {\mathfrak A}$$
where each $S_{j}$ is spanned by the points $p_{\sigma(1)},\dots, p_{\sigma(j+1)}$.
The mapping
$$\sigma \to {\mathcal F}_{\sigma}$$
is a bijection of the set of all permutations of $I_{n}$ to ${\mathfrak A}$.
For any flag ${\mathcal F}\in {\mathfrak A}$
the corresponding permutation of $I_{n}$ will be denoted by $\sigma_{\mathcal F}$;
if
$${\mathcal F}=(p_{i_{1}},S_{1},\dots, S_{n-2},p^{*}_{i_{n+1}})$$
and each $S_{j}$ is the subspace spanned by $p_{i_{1}},\dots, p_{i_{j+1}}$
then $\sigma_{\mathcal F}(j)=i_{j}$ for all $j$,
hence $\sigma^{-1}_{{\mathcal F}}(j)$ is the order number of the minimal
component of the flag ${\mathcal F}$ containing $p_{j}$ and
if this point does not belong to any component of ${\mathcal F}$
then $\sigma^{-1}_{{\mathcal F}}(j)=n+1$.

Denote by ${\mathfrak R}_{ij}$
the set of all flags ${\mathcal F}\in {\mathfrak A}$
satisfying the following condition:

$$1<\sigma^{-1}_{{\mathcal F}}(i)<\sigma^{-1}_{{\mathcal F}}(j)<n+1.$$
If $n=2$ then this set is empty.

\begin{lemma}
The set ${\mathfrak C}_{ij}$ is the disjoint union of
${\mathfrak A}_{i}\cup{\mathfrak A}^{j}$ and ${\mathfrak R}_{ij}$.
\end{lemma}

\begin{proof}
It is trivial that the sets
${\mathfrak A}_{i}\cup{\mathfrak A}^{j}$ and ${\mathfrak R}_{ij}$
are disjoint.
Since $p_{i}$ and $p^{*}_{j}$ do not belong to ${\mathcal X}_{ij}$,
${\mathfrak A}_{i}\cup{\mathfrak A}^{j}$ does not intersects
${\mathfrak X}_{ij}$ and we have
$${\mathfrak A}_{i}\cup{\mathfrak A}^{j}\subset {\mathfrak C}_{ij}.$$
For any flag
$${\mathcal F}=(p_{k},S_{1},\dots,S_{n-2},p^{*}_{m})\in {\mathfrak R}_{ij}$$
the subspace $S_{\sigma^{-1}_{{\mathcal F}}(i)}$
(the minimal subspace of ${\mathcal F}$ containing $p_{i}$)
does not contain the point $p_{j}$;
thus $S_{\sigma^{-1}_{{\mathcal F}}(i)}$ is not an element of ${\mathcal X}_{ij}$ and
${\mathcal F}\notin {\mathfrak X}_{ij}$;
this implies that
$${\mathfrak R}_{ij}\subset {\mathfrak C}_{ij}.$$
Inversely, let us consider a flag ${\mathcal F}\in {\mathfrak C}_{ij}$
which does not belong to ${\mathfrak A}_{i}\cup{\mathfrak A}^{j}$.
Then
$${\mathcal F}\notin {\mathfrak A}_{i}\Rightarrow \sigma_{{\mathcal F}}(1)\ne i
\;\mbox{ and }\;
{\mathcal F}\notin {\mathfrak A}^{j}\Rightarrow \sigma_{{\mathcal F}}(n+1)\ne j.$$
If $\sigma^{-1}_{{\mathcal F}}(i)>\sigma^{-1}_{{\mathcal F}}(j)$
then each component of the flag ${\mathcal F}$ belongs to
${\mathcal X}_{ij}$ and we get ${\mathcal F}\in {\mathfrak X}_{ij}$,
this contradicts to the condition ${\mathcal F}\in{\mathfrak C}_{ij}$.
Therefore, $\sigma^{-1}_{{\mathcal F}}(i)<\sigma^{-1}_{{\mathcal F}}(j)$
and ${\mathcal F}$ is an element of the set ${\mathfrak R}_{ij}$.
\end{proof}

For any subspace $S\in T({\mathfrak A})$ there exists
unique subspace $S^{c}\in T({\mathfrak A})$ non-intersecting $S$ and such that
our projective space is spanned by $S$ and $S^{c}$.
For any flag
$${\mathcal F}=(p_{i}, S_{1},\dots, S_{n-2}, p^{*}_{j})\in {\mathfrak A}$$
the flag
$${\mathcal F}^{c}:=(p_{j}=(p^{*}_{j})^{c},(S_{n-2})^{c},\dots, (S_{1})^{c},
p^{*}_{i}=(p_{i})^{c})\in {\mathfrak A}$$
is said to be the {\it complement}.
The transformation of ${\mathfrak A}$
sending each flag to its complement is bijective;
moreover,
$$({\mathfrak A}_{i})^{c}={\mathfrak A}^{i},\;({\mathfrak A}^{i})^{c}={\mathfrak A}_{i},\;
({\mathfrak R}_{ij})^{c}={\mathfrak R}_{ji}$$
and Lemma 4.4 shows that
\begin{equation}
({\mathfrak C}_{ij})^{c}={\mathfrak C}_{ji}.
\end{equation}

\begin{lemma}
Let $n\ge 3$ and $i,j,k,m$ be distinct elements of $I_{n}$.
Then
$$|{\mathfrak A}^{k}_{m}\cap{\mathfrak R}_{ij}|=\frac{(n-1)!}{2}$$
and
$$|{\mathfrak A}_{k}\cap{\mathfrak R}_{ij}|=|{\mathfrak A}^{k}\cap{\mathfrak R}_{ij}|=
\frac{(n-2)(n-1)!}{2}.$$
\end{lemma}

\begin{proof}
Let us consider the bijective transformation $d_{ij}$ of ${\mathfrak A}$
defined by the formula
$$d_{ij}({\mathcal F}):={\mathcal F}_{\sigma_{ij}\sigma_{\mathcal F}}$$
where $\sigma_{ij}$ is the permutation sending $i$ and $j$ to $j$ and $i$
(respectively) and leaving fixed all elements of the set $I_{n}-\{i,j\}$.
This bijection maps ${\mathfrak A}^{k}_{m}\cap{\mathfrak R}_{ij}$
to ${\mathfrak A}^{k}_{m}-{\mathfrak R}_{ij}$;
thus these sets have the same cardinal number.
Since
$${\mathfrak A}^{k}_{m}=({\mathfrak A}^{k}_{m}\cap{\mathfrak R}_{ij})
\sqcup ({\mathfrak A}^{k}_{m}-{\mathfrak R}_{ij})$$
contains $(n-1)!$ elements,
we get the first equality.
Then the second equality is a consequence of the
following fact: ${\mathfrak A}^{k}\cap{\mathfrak R}_{ij}$
(${\mathfrak A}_{k}\cap{\mathfrak R}_{ij}$)
is the disjoint union of all ${\mathfrak A}^{k}_{l}\cap{\mathfrak R}_{ij}$
(${\mathfrak A}^{l}_{k}\cap{\mathfrak R}_{ij}$)
such that $l\in I_{n}-\{i,j,k\}$.
\end{proof}

\subsection{Dispositions of two complement sets}
Let ${\mathfrak C}_{ij}$ and ${\mathfrak C}_{km}$
be distinct complement sets.
Set
$$J:=\{i,j,k,m\}.$$
Then $2\le |J|\le 4$ and there are the following possibilities for the disposition
of ${\mathfrak C}_{ij}$ and ${\mathfrak C}_{km}$:
\begin{enumerate}
\item[(1)] $i=m$ and $j=k$;
\item[(2)] $|J|=3$ and $i=k$;
\item[(3)] $|J|=3$ and $j=m$;
\item[(4)] $|J|=3$ and $i=m$;
\item[(5)] $|J|=3$ and $j=k$;
\item[(6)]  $|J|=4$.
\end{enumerate}
For each number $l\in \{1,\dots,6\}$ denote
by $n_{l}$ the cardinal number of the intersection of
${\mathfrak C}_{ij}$ and ${\mathfrak C}_{km}$ for the case $(l)$.
If $n=2$ then the case (6) is not realized and the number $n_{6}$
is not defined.

\begin{lemma}
$n_{1}=0$.
\end{lemma}

\begin{proof}
By Lemma 4.4,
$${\mathfrak C}_{ij}=({\mathfrak A}_{i}\cup{\mathfrak A}^{j})\sqcup {\mathfrak R}_{ij},$$
$${\mathfrak C}_{ji}=({\mathfrak A}_{j}\cup{\mathfrak A}^{i})\sqcup {\mathfrak R}_{ji};$$
we have also
$$({\mathfrak A}_{i}\cup{\mathfrak A}^{j})\cap({\mathfrak A}_{j}\cup{\mathfrak A}^{i})
=\emptyset,$$
$$({\mathfrak A}_{i}\cup{\mathfrak A}^{j})\cap {\mathfrak R}_{ji}=
({\mathfrak A}_{j}\cup{\mathfrak A}^{i})\cap {\mathfrak R}_{ij}=\emptyset.$$
and get the equality
$${\mathfrak C}_{ij}\cap {\mathfrak C}_{ji}={\mathfrak R}_{ij}\cap {\mathfrak R}_{ji}.$$
However,
$${\mathcal F}\in {\mathfrak R}_{ij}\Rightarrow
\sigma^{-1}_{{\mathcal F}}(i)<\sigma^{-1}_{{\mathcal F}}(j)$$
contradicts to
$${\mathcal F}\in {\mathfrak R}_{ji}\Rightarrow
\sigma^{-1}_{{\mathcal F}}(j)<\sigma^{-1}_{{\mathcal F}}(i).$$
Therefore, ${\mathfrak R}_{ij}\cap {\mathfrak R}_{ji}=\emptyset$
and the intersection of our complement sets is empty.
\end{proof}

\begin{lemma}
$n_{2}=n_{3}$ and $n_{4}=n_{5}$.
\end{lemma}

\begin{proof}
The equality (4.1) shows that the bijection
$${\mathcal F}\to {\mathcal F}^{c}$$
transfers pairs of complement sets
satisfying the condition (2) or (4) to
pairs of complement sets satisfying the condition (3) or (5), respectively.
\end{proof}

\begin{lemma}
If $n\ge 4$ then
$$n_{2}=n! +(n-2)(n-1)!+\frac{(n-2)(n-3)(n-1)!}{3}.$$
We have also
$$n_{2}=n! +(n-2)(n-1)!$$
if $n=2,3$.
\end{lemma}

\begin{proof}
Since
$${\mathfrak C}_{ij}=({\mathfrak A}_{i}\cup{\mathfrak A}^{j})\sqcup {\mathfrak R}_{ij},$$
$${\mathfrak C}_{im}=({\mathfrak A}_{i}\cup{\mathfrak A}^{m})\sqcup {\mathfrak R}_{im}$$
(Lemma 4.4) and
$$({\mathfrak A}_{i}\cup{\mathfrak A}^{j})\cap
({\mathfrak A}_{i}\cup{\mathfrak A}^{m})={\mathfrak A}_{i},$$
$$({\mathfrak A}_{i}\cup{\mathfrak A}^{j})\cap {\mathfrak R}_{im}=
{\mathfrak A}^{j} \cap {\mathfrak R}_{im},$$
$$({\mathfrak A}_{i}\cup{\mathfrak A}^{m})\cap {\mathfrak R}_{ij}=
{\mathfrak A}^{m}\cap{\mathfrak R}_{ij};$$
the intersection of ${\mathfrak C}_{ij}$ and ${\mathfrak C}_{im}$
is the disjoint union of the following sets:
$${\mathfrak A}_{i},\;{\mathfrak A}^{j}\cap{\mathfrak R}_{im},\;
{\mathfrak A}^{m}\cap{\mathfrak R}_{ij},\;
{\mathfrak R}_{ij}\cap {\mathfrak R}_{im}.$$
If $n=2$ then the sets ${\mathfrak R}_{ij}$ and ${\mathfrak R}_{im}$
are empty and the required equality is trivial.
For the case when $n\ge 3$ we have
$$n_{2}=|{\mathfrak A}_{i}|+|{\mathfrak A}^{j}\cap{\mathfrak R}_{im}|+
|{\mathfrak A}^{m}\cap{\mathfrak R}_{ij}|+|{\mathfrak R}_{ij}\cap {\mathfrak R}_{im}|
=n! +(n-2)(n-1)!+|{\mathfrak R}_{ij}\cap {\mathfrak R}_{im}|$$
(Lemma 4.5).
We will find the cardinal number of the set ${\mathfrak R}_{ij}\cap {\mathfrak R}_{im}$.

First of all note that this set is empty if $n=3$,
since for this case
for any flag ${\mathcal F}\in {\mathfrak A}$
one of the numbers
\begin{equation}
\sigma^{-1}_{{\mathcal F}}(i),\;\sigma^{-1}_{{\mathcal F}}(j),\;
\sigma^{-1}_{{\mathcal F}}(m)
\end{equation}
is equal to $1$ or $n+1$.
Let $n\ge 4$.
Then $I_{n}-J$ contains at list two elements and
we take distinct numbers $l$ and $s$ belonging to this set.
If ${\mathcal F}\in {\mathfrak A}^{s}_{l}$ then
there are exact $6$ possibilities for the disposition of
the numbers (4.2) in the set $\sigma^{-1}_{{\mathcal F}}(J)$;
only for two of these cases $\sigma^{-1}_{{\mathcal F}}(i)$ is less
than both $\sigma^{-1}_{{\mathcal F}}(j)$, $\sigma^{-1}_{{\mathcal F}}(m)$
and the flag ${\mathcal F}$ belongs to ${\mathfrak R}_{ij}\cap {\mathfrak R}_{im}$.
Therefore,
$$|{\mathfrak A}^{s}_{l}\cap{\mathfrak R}_{ij}\cap {\mathfrak R}_{im}|=
\frac{2|{\mathfrak A}^{s}_{l}|}{6}=\frac{(n-1)!}{3}$$
The set ${\mathfrak R}_{ij}\cap {\mathfrak R}_{im}$
is the disjoint union of all
${\mathfrak A}^{s}_{l}\cap{\mathfrak R}_{ij}\cap {\mathfrak R}_{im}$
such that $l$ and $s$ do not belong to $J$.
It is easy to see that there are $(n-2)(n-3)$ distinct
${\mathfrak A}^{s}_{l}$ with $s,l\in I_{n}-J$.
Thus
$$|{\mathfrak R}_{ij}\cap {\mathfrak R}_{im}|=\frac{(n-2)(n-3)(n-1)!}{3}.$$
\end{proof}

\begin{lemma}
If $n\ge 4$ then
$$n_{4}=(n-1)! +(n-2)(n-1)!+\frac{(n-2)(n-3)(n-1)!}{6}$$
For the case when $n=2$ or $3$ we have
$$n_{4}=(n-1)! +(n-2)(n-1)!.$$
\end{lemma}

\begin{proof}
By Lemma 4.4,
$${\mathfrak C}_{ij}=({\mathfrak A}_{i}\cup{\mathfrak A}^{j})\sqcup {\mathfrak R}_{ij},$$
$${\mathfrak C}_{ki}=({\mathfrak A}_{k}\cup{\mathfrak A}^{i})\sqcup {\mathfrak R}_{ki}.$$
A direct verification shows that
$$({\mathfrak A}_{i}\cup{\mathfrak A}^{j})\cap({\mathfrak A}_{k}\cup{\mathfrak A}^{i})=
{\mathfrak A}^{j}_{k},$$
$$({\mathfrak A}_{i}\cup{\mathfrak A}^{j})\cap {\mathfrak R}_{ki}=
{\mathfrak A}^{j}\cap{\mathfrak R}_{ki},$$
$$({\mathfrak A}_{k}\cup{\mathfrak A}^{i})\cap {\mathfrak R}_{ij}=
{\mathfrak A}_{k}\cap {\mathfrak R}_{ij}.$$
Hence ${\mathfrak C}_{ij}\cap {\mathfrak C}_{ki}$ is the disjoint union of
the sets
$${\mathfrak A}^{j}_{k},\;{\mathfrak A}^{j}\cap{\mathfrak R}_{ki},\;
{\mathfrak A}_{k}\cap {\mathfrak R}_{ij},\; {\mathfrak R}_{ij}\cap{\mathfrak R}_{ki}$$
and for the case when $n=2$ the required equality is trivial
(since the sets ${\mathfrak R}_{ij}$ and ${\mathfrak R}_{ki}$ are empty).
If $n\ge 3$ then
$$n_{4}=|{\mathfrak A}^{j}_{k}|+|{\mathfrak A}^{j}\cap{\mathfrak R}_{ki}|+
|{\mathfrak A}_{k}\cap {\mathfrak R}_{ij}|+|{\mathfrak R}_{ij}\cap {\mathfrak R}_{ki}|
=(n-1)! +(n-2)(n-1)!+|{\mathfrak R}_{ij}\cap {\mathfrak R}_{ki}|.$$
The set ${\mathfrak R}_{ij}\cap {\mathfrak R}_{ki}$ is empty if $n\le 3$;
indeed, for any flag ${\mathcal F}\in {\mathfrak A}$
one of the numbers
\begin{equation}
\sigma^{-1}_{{\mathcal F}}(i),\;\sigma^{-1}_{{\mathcal F}}(j),\;
\sigma^{-1}_{{\mathcal F}}(k)
\end{equation}
is equal to $1$ or $n+1$.
Assume that $n\ge 4$ (for this case $I_{n}-J$ contains at least two elements) and
take two distinct numbers $l$ and $s$ belonging to $I_{n}-J$.
It was noted above (the proof of Lemma 4.8)
that for any flag ${\mathcal F}\in {\mathfrak A}^{s}_{l}$
there are exact $6$ possibilities for the disposition of
the numbers (4.3) in the set $\sigma^{-1}_{{\mathcal F}}(J)$.
The inequality
$$\sigma^{-1}_{{\mathcal F}}(k)<\sigma^{-1}_{{\mathcal F}}(i)<\sigma^{-1}_{{\mathcal F}}(j)$$
holds only for one of these cases.
Then we have
$$|{\mathfrak A}^{s}_{l}\cap{\mathfrak R}_{ij}\cap {\mathfrak R}_{ki}|=
\frac{|{\mathfrak A}^{s}_{l}|}{6}=\frac{(n-1)!}{6}$$
and the arguments used to prove Lemma 4.8 show that
$$|{\mathfrak R}_{ij}\cap {\mathfrak R}_{ki}|=\frac{(n-2)(n-3)(n-1)!}{6}.$$
\end{proof}

\begin{lemma}
If $n\ge 5$ then
$$n_{6}=n!+ (n-1)(n-1)!+\frac{(n-3)(n-4)(n-1)!}{4}.$$
We have also
$$n_{6}=n!+ (n-1)(n-1)!$$
if $n=3,4$.
\end{lemma}

Recall that for the case when $n=2$ the number $n_{6}$ is not defined.

\begin{proof}
Lemma 4.4 says that
$${\mathfrak C}_{ij}=({\mathfrak A}_{i}\cup{\mathfrak A}^{j})\sqcup {\mathfrak R}_{ij},$$
$${\mathfrak C}_{km}=({\mathfrak A}_{k}\cup{\mathfrak A}^{m})\sqcup {\mathfrak R}_{km}.$$
Since
$$({\mathfrak A}_{i}\cup{\mathfrak A}^{j})\cap({\mathfrak A}_{k}\cup{\mathfrak A}^{m})=
{\mathfrak A}^{m}_{i}\sqcup {\mathfrak A}^{j}_{k},$$
${\mathfrak C}_{ij}\cap{\mathfrak C}_{km}$ is
the union of the following non-intersecting sets:
$${\mathfrak A}^{m}_{i},\;{\mathfrak A}^{j}_{k},\;
({\mathfrak A}_{i}\cup{\mathfrak A}^{j})\cap {\mathfrak R}_{km},\;
({\mathfrak A}_{k}\cup{\mathfrak A}^{m})\cap {\mathfrak R}_{ij},\;
{\mathfrak R}_{ij}\cap {\mathfrak R}_{km}.$$
The set ${\mathfrak A}_{i}\cup{\mathfrak A}^{j}$ is the disjoint union of all
${\mathfrak A}^{l}_{i}$, $l\ne i$ and all ${\mathfrak A}^{j}_{s}$, $s\ne j$.
Each of these collections contains $n$ sets,
but the set ${\mathfrak A}^{j}_{i}$ belongs to the both collections;
note also that ${\mathfrak A}^{m}_{i}$ and ${\mathfrak A}^{j}_{k}$ do not
intersect ${\mathfrak R}_{km}$.
In other words, $({\mathfrak A}_{i}\cup{\mathfrak A}^{j})\cap {\mathfrak R}_{km}$
is the disjoint union of $2n-3$ distinct sets containing $\frac{(n-1)!}{2}$ elements
(Lemma 4.5) and
$$|({\mathfrak A}_{i}\cup{\mathfrak A}^{j})\cap {\mathfrak R}_{km}|=
\frac{(2n-3)(n-1)!}{2}.$$
Similarly,
$$|({\mathfrak A}_{k}\cup{\mathfrak A}^{m})\cap {\mathfrak R}_{ij}|=
\frac{(2n-3)(n-1)!}{2}.$$
Therefore,
$$n_{6}=|{\mathfrak A}^{m}_{i}|+|{\mathfrak A}^{j}_{k}|+
|({\mathfrak A}_{i}\cup{\mathfrak A}^{j})\cap {\mathfrak R}_{km}|+
|({\mathfrak A}_{k}\cup{\mathfrak A}^{m})\cap {\mathfrak R}_{ij}|+
|{\mathfrak R}_{ij}\cap {\mathfrak R}_{km}|=$$
$$2(n-1)!+(2n-3)(n-1)!+|{\mathfrak R}_{ij}\cap {\mathfrak R}_{km}|=
n!+(n-1)(n-1)!+|{\mathfrak R}_{ij}\cap {\mathfrak R}_{km}|.$$
Thus we have to find the cardinal number of the set
${\mathfrak R}_{ij}\cap {\mathfrak R}_{km}$.

If $n\le 4$ then ${\mathfrak R}_{ij}\cap {\mathfrak R}_{km}$ is empty,
since for this case for any flag ${\mathcal F}\in {\mathfrak A}$
at least one of the numbers
$$\sigma^{-1}_{{\mathcal F}}(i),\;\sigma^{-1}_{{\mathcal F}}(j),\;
\sigma^{-1}_{{\mathcal F}}(k),\;\sigma^{-1}_{{\mathcal F}}(m)$$
is equal to $1$ or $n+1$.
Now assume that $n\ge 5$.
Then $I_{n}-J$ contains at least two elements and
we take two distinct numbers $l$ and $s$ belonging to this set.
The bijection $d_{ij}$ (see the proof of Lemma 4.5) transfers
$${\mathfrak A}^{s}_{l}\cap{\mathfrak R}_{km}\cap{\mathfrak R}_{ij}\;
\mbox{ to }\;({\mathfrak A}^{s}_{l}\cap{\mathfrak R}_{km})-{\mathfrak R}_{ij}$$
Hence these sets have the same cardinal number.
Since
$${\mathfrak A}^{s}_{l}\cap{\mathfrak R}_{km}=
({\mathfrak A}^{s}_{l}\cap{\mathfrak R}_{km}\cap{\mathfrak R}_{ij})\sqcup
(({\mathfrak A}^{s}_{l}\cap{\mathfrak R}_{km})-{\mathfrak R}_{ij})$$
contains $\frac{(n-1)!}{2}$ elements, we have
$$|{\mathfrak A}^{s}_{l}\cap {\mathfrak R}_{ij}\cap {\mathfrak R}_{km}|=
\frac{(n-1)!}{4}.$$
The set ${\mathfrak R}_{ij}\cap {\mathfrak R}_{km}$
is the disjoint unions of all
${\mathfrak A}^{s}_{l}\cap {\mathfrak R}_{ij}\cap {\mathfrak R}_{km}$
such that $s,l\in I_{n}-J$.
There are exactly $(n-3)(n-4)$ distinct sets satisfying this condition
and we get the equality
$$|{\mathfrak R}_{ij}\cap {\mathfrak R}_{km}|=
\frac{(n-3)(n-4)(n-1)!}{4}.$$
\end{proof}

\begin{lemma}
The number $n_{2}=n_{3}$ is not equal to $n_{1}$, $n_{4}=n_{5}$ and $n_{6}$.
\end{lemma}

\begin{proof}
Clearly, $n_{2}> 0=n_{1}$.
It is easy to see that $n_{2}>n_{4}$.
We need to establish that $n_{2}$ is not equal to $n_{6}$ if $n\ge 3$.

If $n=3$ or $4$ then
$$n_{6}-n_{2}=(n-1)!=2$$
or
$$n_{6}-n_{2}=(n-1)!-\frac{(n-2)(n-3)(n-1)!}{3}=2,$$
respectively.
For the case when $n\ge 5$ we have
$$n_{2}-n_{6}=\frac{(n-2)(n-3)(n-1)!}{3}-(n-1)!-\frac{(n-3)(n-4)(n-1)!}{4}=$$
$$(n-1)!\frac{4(n-2)(n-3)-3(n-3)(n-4)-12}{12}=(n-1)!\frac{n^{2}+n-24}{12}>0.$$
\end{proof}

We will say that the complement sets ${\mathfrak C}_{ij}$ and ${\mathfrak C}_{km}$
are {\it adjacent} if one of the cases (2) or (3) is realized, i.e. $i=k$ or $j=m$.

\begin{lemma}
The mapping $f$ transfers complement subsets of ${\mathfrak A}$ to
complement subsets of ${\mathfrak A}'$; moreover,
adjacent complement subsets go over to
adjacent complement subsets.
\end{lemma}

\begin{proof}
The first statement follows from Lemma 4.3.
Since $f$ is injective, the $f$-images of two adjacent complement subsets of
${\mathfrak A}$ are complement subsets of ${\mathfrak A}'$
which intersection contains $n_{2}=n_{3}$ elements.
By Lemma 4.11, these complement subsets are adjacent.
\end{proof}

\subsection{Proof of Lemma 4.1}
First of all we give two simple lemmas.

\begin{lemma}
For each number $i\in I_{n}$ we have
$$\bigcap_{j\in I_{n}-\{i\}}{\mathfrak C}_{ij}={\mathfrak A}_{i}$$
and
$$\bigcap_{j\in I_{n}-\{i\}}{\mathfrak C}_{ji}={\mathfrak A}^{i}.$$
\end{lemma}

\begin{proof}
It is trivial that the intersection of all ${\mathfrak C}_{ij}$, $j\in I_{n}-\{i\}$
contains ${\mathfrak A}_{i}$.
If ${\mathcal F}\in {\mathfrak A}_{k}$ and $k\ne i$ then
${\mathcal F}$ does not belong to ${\mathfrak C}_{ik}$.
Therefore, our intersection coincides with ${\mathfrak A}_{i}$
and we get the first equality.
The transformation
$${\mathcal F} \to {\mathcal F}^{c}$$
sends the first equality to the second equality.
\end{proof}

\begin{lemma}
For any collection of $n$ mutually adjacent complement subsets
of ${\mathfrak A}$ there exists a number $i\in I_{n}$ such that
this collection is consisting of all ${\mathfrak C}_{ij}$ or all
${\mathfrak C}_{ji}$ with $j\in I_{n}-\{i\}$.
\end{lemma}

\begin{proof}
Easy verification.
\end{proof}

Now we can prove Lemma 4.1.

By Lemma 4.13,
$$f({\mathfrak A}_{i})=\bigcap_{j\in I_{n}-\{i\}}f({\mathfrak C}_{ij})$$
for any $i\in I_{n}$.
Lemma 4.12 shows that $\{f({\mathfrak C}_{ij})\}_{j\in I_{n}-\{i\}}$
is a collection of $n$ mutually adjacent complement subsets of ${\mathfrak A}'$.
Then Lemma 4.13 and 4.14 guarantee the existence a number
$\sigma(i)\in I_{n}$ such that $f({\mathfrak A}_{i})$ coincides with
${\mathfrak A}'_{\sigma(i)}$ or ${\mathfrak A}'^{\sigma(i)}$.

We want to establish that the set
$$J:=\{\;i\in I_{n}\;|\;f({\mathfrak A}_{i})={\mathfrak A}'_{\sigma(i)}\;\}$$
is empty or coincides with $I_{n}$.
In particular, this means that the mapping $\sigma:I_{n}\to I_{n}$ is bijective
(since the restriction of $f$ to ${\mathfrak A}$ is a bijection to ${\mathfrak A}'$).

\begin{proof}
Assume that $J$ is not empty and does not coincide with $I_{n}$.
There exist $k\in J$ and $m\in I_{n}-J$ such that $\sigma(k)\ne \sigma(m)$
(otherwise, the mapping $\sigma:I_{n}\to I_{n}$ is constant and $f({\mathfrak A})$
is contained in certain ${\mathfrak A}'_{l}\cup {\mathfrak A}'^{l}$).
Then
$$f({\mathfrak A}_{k}\cap {\mathfrak A}_{m})=f({\mathfrak A}_{k})\cap f({\mathfrak A}_{m})=
{\mathfrak A}'_{\sigma(k)}\cap {\mathfrak A}'^{\sigma(m)}=
{\mathfrak A}'^{\sigma(m)}_{\sigma(k)}\ne\emptyset$$
contradicts to ${\mathfrak A}_{k}\cap {\mathfrak A}_{m}=\emptyset$.
\end{proof}

If $J$ coincides with $I_{n}$ then
$f({\mathfrak A}^{i})={\mathfrak A}'^{\sigma(i)}$
for each $i\in I_{n}$.

\begin{proof}
The arguments given above imply the existence of a permutation
$\omega$ of $I_{n}$ such that
$$
f({\mathfrak A}^{i})={\mathfrak A}'_{\omega(i)}\;\;\;\;\;\;
\forall\;i\in I_{n}
$$
or
$$
f({\mathfrak A}^{i})={\mathfrak A}'^{\omega(i)}\;\;\;\;\;\;
\forall\;i\in I_{n}.
$$
Clearly, for the first case the mapping $f$ is not injective.
Thus the second equality holds true.
Since ${\mathfrak A}_{i}\cap {\mathfrak A}^{i}$ is empty,
$$f({\mathfrak A}_{i}\cap {\mathfrak A}^{i})=f({\mathfrak A}_{i})\cap f({\mathfrak A}^{i})=
{\mathfrak A}'_{\sigma(i)}\cap{\mathfrak A}'^{\omega(i)}=\emptyset ;$$
this means that $\sigma(i)=\omega(i)$.
\end{proof}

If $J$ is empty then the similar arguments show that
$f({\mathfrak A}^{i})={\mathfrak A}'_{\sigma(i)}$
for all $i\in I_{n}$.

\subsection{The final part of the proof of Theorem 3.1}
Recall that (Example 3.1) for any subspaces $S\subset P$ and $S'\subset P'$
we denote by ${\mathfrak F}(S)$ and ${\mathfrak F}'(S')$
the sets of all maximal flags containing these subspaces.

Let $p\in P$. Assume that
${\mathcal F}_{1}$ and ${\mathcal F}_{2}$ are distinct flags belonging to the set
${\mathfrak F}(p)$.
Since apartments containing ${\mathcal F}_{1}$ and ${\mathcal F}_{2}$
exist, Lemma 4.1 shows that
the flags $f({\mathcal F}_{1})$ and $f({\mathcal F}_{2})$
have the same point or the same $(n-1)$-dimensional subspace.

Now take ${\mathcal F}_{0}\in{\mathfrak F}(p)$
and denote by ${\mathfrak X}$ (or ${\mathfrak Y}$)
the set of all flags ${\mathcal F}\in{\mathfrak F}(p)$
such that $f({\mathcal F}_{0})$ and $f({\mathcal F})$ have
the same point (or the same $(n-1)$-dimensional subspace, respectively).
Then
\begin{equation}
{\mathfrak X}\cup {\mathfrak Y}={\mathfrak F}(p).
\end{equation}
We want to show that one of these sets coincides with ${\mathfrak F}(p)$.

\begin{proof}
If
\begin{equation}
{\mathfrak X}\ne{\mathfrak F}(p)\;\mbox{ and }\;{\mathfrak Y}\ne{\mathfrak F}(p)
\end{equation}
then (4.4) implies the existence of
$${\mathcal F}_{1}\in{\mathfrak X}-{\mathfrak Y}\;\mbox{ and }\;
{\mathcal F}_{2}\in{\mathfrak Y}-{\mathfrak X}.$$
The points and the $(n-1)$-dimensional subspaces
contained in the flags $f({\mathcal F}_{1})$ and $f({\mathcal F}_{2})$ are different.
Hence one of the conditions (4.5) does not hold.
\end{proof}

Thus we have established that for any point $p\in P$ there exists a point
$g(p)$ belonging to $P'$ or $P'^{*}$ and such that
$$f({\mathfrak F}(p))\subset{\mathfrak F}'(g(p)).$$
Since any two points of $P$ are contained in some base,
Lemma 4.1 guarantees that one of the following two possibilities
is realized:
\begin{enumerate}
\item[(A)] $g(p)\in P'$ for all $p\in P$,
\item[(B)] $g(p)\in P'^{*}$ for all $p\in P$.
\end{enumerate}
Similarly, for any point $p^{*}\in P^{*}$
there is a point $h(p^{*})$ of $P'$ or $P'^{*}$ such that
$$f({\mathfrak F}(p^{*}))\subset{\mathfrak F}'(h(p^{*})).$$
Moreover, it follows from Lemma 4.1 that for the cases (A) and (B) we have
$$h(p^{*})\in P'^{*}\;\;\;\forall\;p^{*}\in P^{*}\;
\mbox{ and }\;
h(p^{*})\in P'\;\;\;\forall\;p^{*}\in P^{*},$$
respectively.

\begin{lemma}
If $p\in P$ and $p^{*}\in P^{*}$ contains $p$
then $g(p)$ and $h(p^{*})$ are incident subspaces
\footnote{recall that two subspaces are incident if one of them is contained to
the other}.
\end{lemma}

\begin{proof}
Consider a maximal flag ${\mathfrak F}$ containing $p$ and $p^{*}$.
Then $g(p)$ and $h(p^{*})$ are contained in the flag $f({\mathfrak F})$.
This gives the required.
\end{proof}

The mapping $g$ is a strong embedding of $P$ to $P'$ or $P'^{*}$.

\begin{proof}
It is trivial that $g$ maps bases to bases.
Hence it transfers independent subsets to independent subsets
(a subset of a projective space is independent if and only if it is contained
in certain base for this space).
In particular, this means that $g$ is injective.
Thus we have to show that $g$ is collinearity preserving.
We prove this statement for the case (A),
the case (B) is similar.

Any two distinct points $p_{1},p_{2}\in P$
are contained in some base $B=\{p_{i}\}^{n+1}_{i=1}$ for $P$.
The points
$$p'_{1}:=g(p_{1}),\dots, p'_{n+1}:=g(p_{n+1})$$
form a base for $P'$;
in what follows this base will be denoted by $B'$.
For each $i\in I_{n}$ we set
$$p^{*}_{i}:=\overline{B-\{p_{i}\}},\;\;p'^{*}_{i}:=\overline{B'-\{p'_{i}\}}$$
and denote by ${\mathfrak A}$ and ${\mathfrak A}'$ the apartments in ${\mathfrak F}$
and ${\mathfrak F}'$ associated with the bases $B$ and $B'$, respectively.
Then $f({\mathfrak A})={\mathfrak A}'$ and Lemma 4.1 shows that
$h(p^{*}_{i})=p'^{*}_{i}$ for all $i\in I_{n}$.
Since
$$p_{1}p_{2}=\bigcap^{n+1}_{i=3}p^{*}_{i},$$
we have
$$g(p_{1}p_{2})\subset \bigcap^{n+1}_{i=3}h(p^{*}_{i})=
\bigcap^{n+1}_{i=3}p'^{*}_{i}=p'_{1}p'_{2}$$
by Lemma 4.15.
\end{proof}

Lemma 4.15 states that the $g$-image of any $p^{*}\in P^{*}$
is contained in $h(p^{*})$.
Since $g$ is a strong embedding to $P'$ or $P'^{*}$,
$\overline{g(p^{*})}$ is an $(n-1)$-dimensional subspaces of $P'$ or $P'^{*}$
(respectively). The similar holds for $h(p^{*})$.
Hence the subspaces $h(p^{*})$ and $\overline{g(p^{*})}$ are coincident.
This means that $h=g^{*}$ (Subsection 2.3).
In particular, Theorem 3.1 is proved for the case when $n=2$.

Now assume that $n\ge 3$ and prove Theorem 3.1 by induction.
Take an arbitrary point $p\in P$ and
consider the restriction of $f$ to ${\mathfrak F}(p)$.
It is an injection to ${\mathfrak F}'(g(p))$
sending apartments to apartments (Example 3.1).
By the inductive hypothesis, this mapping is induced
by a strong embedding $t$ of ${\mathcal L}_{p}$ to ${\mathcal L}'_{g(p)}$ or
${\mathcal L}'^{*}_{g(p)}$ (Example 2.5).
It is easy to see that the dual embedding $t^{*}$ coincides
with the restriction of $h$ to $P^{*}_{p}$; denote this restriction by $h'$.
Then $t=h'^{*}$ is defined by $g$ (Example 2.5).
This implies the required: the mapping $f$ is induced by the embedding $g$.

{\bf Acknowledgment}. The author thanks H.~Havlicek for his interest and valuable discussions .

\end{document}